\documentclass[12pt]{amsart}
\usepackage{amssymb,amsfonts,amscd}
\theoremstyle{plain}
\newtheorem{thm}{Theorem}[section]
\newtheorem{cor}[thm]{Corollary}

\newtheorem{lem}[thm]{Lemma}
\newtheorem{property}[thm]{Condition}

\theoremstyle{definition}
\newtheorem*{defn}{Definition}

\numberwithin{equation}{section}

\let\le\leqslant

\makeatletter
\renewcommand\sum{\DOTSB\sum@\limits}
\renewcommand\prod{\DOTSB\prod@\limits}
\makeatother
\newcommand\blank{\mathord{\hbox to 1.5ex{\hrulefill}}\,}
\let\opn\operatorname

\let\mf\mathfrak
\def\q{\mf q}
\newcommand\Z{\mathbb Z}
\newcommand\N{\mathbb N}
\let\ph\varphi
\let\eps\varepsilon

\title[A new look at the decomposition of unipotents]
{A new look at the decomposition of unipotents and the normal structure of Chevalley groups}

\author{Alexei Stepanov}

\address{St. Petersburg State University, Department of Mathematics and Mechanics}
\address{St. Petersburg State Electrotechnical University ``LETI''}

\email{stepanov239@gmail.com}

\thanks
{The work on this publication was supported by
Russian Science Foundation, grant no.~14-11-00297}



\begin{document}
\begin{abstract}
The current article continues a series of papers on the decomposition of unipotents in
a Chevalley group $\opn{G}(\Phi,R)$ over a commutative ring $R$ with a reduced irreducible
root system $\Phi$.
Fix $h\in\opn{G}(\Phi,R)$. Let us call an element $a\in\opn{G}(\Phi,R)$ ``good'',
if it belongs to the unipotent radical of a parabolic subgroup and the conjugate to $a$ by $h$
lies in another parabolic subgroup (all parabolic subgroups are assumed to contain the same
split maximal torus).
The ``decomposition of unipotents'' method is a representation of an elementary root unipotent
element as a product of ``good'' elements.
Decomposition of unipotents implies a simple proof of normality of the elementary subgroup
and the standardness of the normal structure of $\opn{G}(\Phi,R)$. However, such a decomposition
is available not for all root systems. In this article we show that to prove the
standardness of the normal structure it suffices to find one ``good'' element
for the generic element of the group scheme $\opn{G}(\Phi,\blank)$. We also construct some
``good'' elements. The question: ``When good elements span the elementary subgroup?''
will be considered in a subsequent article of the series.
\end{abstract}

\maketitle

\section*{Introduction}

The goal of the current article is to give a new proof of the normal structure theorem for
a Chevalley group provided that the structure constants are invertible.
Let $G=\opn{G}(\Phi,\blank)$ be a Chevalley--Demazure group scheme with a reduced
irreducible root system $\Phi$ and $E=\opn{E}(\Phi,\blank)$ its elementary subgroup.
Let $R$ be a commutative ring with a unit. For an ideal $\mf a$ of $R$ denote by
$C(R,\mf a)$ the full congruence subgroup of level $\mf a$ and let
$E(R,\mf a)$ denote the relative elementary subgroup (see the next section for precise definitions
and notation).

\begin{defn}
We say that the normal structure of the group $G(R)$ is standard if for each subgroup $H\le G(R)$,
normalized by $E(R)$, there exists a unique ideal $\mf a$ of $R$ such that
$$E(R,\mf a)\le H\le C(R,\mf a).$$
\end{defn}

Our proof of the standardness of the normal structure of $G(R)$
is closely related to the decomposition of unipotents techniques.
Fix a split maximal torus $T$ in $G$. We always assume that all root and parabolic
subschemes correspond to $T$.
Let $h\in G(R)$, $r\in R$, and $\alpha\in\Phi$.
Choose a parabolic subscheme $P$ such that the root subscheme
$X_\alpha$ is contained in the unipotent radical $U_P$ of $P$.
Suppose that there exist elements $a_1,\dots,a_m\in U_P(R)$ whose product is equal to
$x_\alpha(r)$ and such that $a_i^h\in Q_i(R)$ for all $i=1,\dots,m$,
where $Q_i$ are parabolic subschemes and $a_i^h=h^{-1}a_ih$.
The representation
$$
x_\alpha(r)^h=a_1^h\dots a_m^h\in Q_1(R)\dots Q_m(R)
$$
is called the decomposition of the unipotent.

Usually, a variation of Whitehead--Vaserstein lemma allows to decompose elements
$a_i^h$ into a product of elementary root unipotents, see e.\,g.~\cite{VasersteinStability,StepVavDecomp}.
This provides a straightforward proof of the normality of the elementary subgroup along
with an efficient bound on the length of $x_\alpha(r)^h$ in elementary generators,
which is much better than a bound that can be obtained by a localization procedure.

A scheme of the normal structure theorem via decomposition of unipotents is as follows,
see~\cite{StepVavDecomp}. Using the standard commutator formulas one reduces the problem to
extraction an elementary root unipotent from a noncentral subgroup $H\le G(R)$,
normalized by $E(R)$ (see e.\,g.~\cite{VavStav}). If $h\in H$ is a noncentral element and
$a\in U_P(R)$ is such that $a^h$ belongs to a parabolic subgroup,
then it is possible to produce an elementary root unipotent out of the commutator $[a,h]\in H$.
This unipotent is nontrivial as soon as $[a,h]\ne e$.
Decomposition of unipotents ensures that elements $a$ satisfying the above property span the
elementary subgroup. Therefore, a noncentral element $h$ cannot commute with all such elements $a$
and we are done.

Usually it is not very difficult to find an element $a$ satisfying the condition
indicated above. The real challenge is to prove that such elements generate the
elementary subgroup and this is not always true. Therefore, it is important to establish
the normal structure theorem based only on the existence of the ``good'' element.

The idea of the proof is inspired by the generic element method, which was named the
``universal localization method'' in~\cite{StepUniloc}. Clearly, instead of an arbitrary
element $h\in G(R)$, where $R$ is an arbitrary ring, it suffices to consider the generic
element $g\in G(A)$ of the scheme $G$ over the affine algebra $A=\Z[G]$.
Suppose that we managed to find a nontrivial element $a\in U_P(A)$ such that $a^g$ lies in
a parabolic subgroup $Q(A)$. We have already noticed that a nontrivial root unipotent can be
produced out of a noncentral element from a parabolic subgroup.
Since the generic element commutes only with the central elements,
any subgroup that contains $g$ and is normalized by $E(A)$
contains a root element $x_\beta(t)$ for some $\beta\in\Phi$ è $t\in A\smallsetminus\{0\}$.

Now, let  $R$ be a ring and let $H$ be a subgroup of $G(R)$ normalized by $E(R)$.
Sending $g$ to an arbitrary element of $H$ one notices that the image of $x_\beta(t)$
under this homomorphism lies in $H$.
We get an alternative: either $H$ contains $x_\beta(r)$ for some $r\in R\smallsetminus\{0\}$
and we are done, or $H$ is contained in a subscheme $S$ defined by the equation $t=0$.
The construction of $a$ ensures that $S(F)\ne G(F)$ for any field $F$.
Therefore the image of $H$ under the map induced by a homomorphism
$R\to F$ is a proper normal subgroup and, therefore, lies in the center.
Then, the centrality of $H$ follows by the standard radical reduction.

Note that the proof of the normal structure theorem from the current article
does not depend on the result over a field or a local ring.
The case of a ground field as well as the radical reduction follows immediately from
the techniques of producing unipotents, developed in the article.

\section{Notation}
Let $a,b,c$ be elements of a group $G$. Denote by $a^b=b^{-1}ab$ the element conjugate to
$a$ by $b$. The commutator $a^{-1}b^{-1}ab$ is denoted by $[a,b]$.
Let $S$ be a subset of $G$. By $\langle S\rangle$ we denote the subgroup
spanned by $S$. For  subgroups $A$ and $B$ of $G$ by $A^B$ we denote the subgroup of $G$
generated by $a^b$ for all $a\in A$ and $b\in B$. In other words, $A^B$ is the smallest subgroup
containing $A$ and normalized by $B$.
The mutual commutator subgroup $[A,B]$ is a subgroup of $G$, generated by all the commutators
$[a,b]$, $a\in A$, $b\in B$.

All rings and algebras are assumed to be commutative and to contain a unit.
All homomorphisms preserve unit elements.
The multiplicative group of $R$ is denoted by $R^*$.
Let $s\in R$. The principal localization at the element $s$
(i.\,e. the localization at the multiplicative subset generated by $s$)
is denoted by $R_s$.

Let $K$ be a ring and $G$ an affine group scheme over $K$.
Denote by $A=K[G]$ the affine algebra of the scheme $G$.
By the definition of an affine scheme, an element $h\in G(R)$ can be identified with a
homomorphism $h:A\to R$. We always do this identification, i.\,e. we always view elements
of the group of points $G(R)$ of the scheme $G$ over a $K$-algebra $R$ as homomorphisms
from $A$ to $R$. Denote by $g\in G(A)$ the generic element of the scheme $G$, i.\,e.
the identity map $\operatorname{id}_A:A\to A$. An element $h\in G(R)$ induces the homomorphism
$G(h):G(A)\to G(R)$ by the rule $G(h)(a)=h\circ a$ for all $a\in G(A)$.
It follows that the image of $g$ under the action of $G(h)$ is equal to~$h$. In the rest of the article
for a ring homomorphism $\ph:R\to R'$ the induced group homomorphism
$G(\ph):G(R)\to G(R')$ is denoted again by $\ph$.
This cannot lead to a confusion as one always can determine the meaning of $\ph$
by the argument type of this homomorphism.
In view of this agreement we have $h(g)=h\circ\operatorname{id}_A=h$.
If $R$ is a subring of $R'$, then we usually identify elements of $G(R)$ with their canonical
images in $G(R')$.
\emph{The notation $A$ and $g$ introduced above is kept throughout the article}.

In what follows $G=\opn{G}(\Phi,\blank)$ denotes a Chevalley--Demazure group scheme over $\Z$
with a reduced irreducible root system $\Phi$ of rank at least 2.
By definition there exists a split maximal torus $T$ in $G$.
Fix $T$; by default all parabolic subgroups contain $T$ and all root subgroups are normalized
by $T$.
Let $E$ be the elementary subgroup of $G$, i.\,e. the span of all root subgroups.

The center of the group $G(R)$ is denoted by $C(R)$.
Let $\mf a$ be an ideal of a ring $R$. As usually, the principal congruence subgroup $G(R,\mf a)$ is the
kernel of the reduction homomorphism $\rho_{\mf a}:G(R)\to G(R/\mf a)$,
whereas the full congruence subgroup $C(R,\mf a)$ is the inverse image of
$C(R/\mf a)$ under this homomorphism.
The relative elementary group $E(R,\mf a)$ is the normal closure in $E(R)$ of the
subgroup, generated by all root unipotents $x_\alpha(r)$, $\alpha\in\Phi$, $r\in\mf a$.

The main result of the current article is proved under the following condition.
\begin{property}\label{StrConstants}
If $\Phi$ is doubly laced (i.\,e. $\Phi=B_l,C_l$ or $F_4$), then $2\in R^*$.
If $\Phi=G_2$, then $3\in R^*$ and $R$ has no residue fields of two elements.
\end{property}

\section{Bruhat and Gauss decompositions}\label{Gauss}
Fix a Borel subgroup $B$ containing the torus $T$.
Let $U$ and $U^-$ denote the unipotent radicals of $B$ and the opposite Borel subgroup $B^-$
respectively. Denote by $W$ the Weyl group of the root system $\Phi$.
Let $N=N_G(T)$ be the scheme theoretic normalizer of the torus $T$, i.\,e. the largest subscheme in $G$
such that $N(R)$ normalizes $T(R)$ for any ring $R$
(by~\cite[Exp.~11, Corollaire~5.3~bis]{SGA} such a subscheme exists).
The quotient scheme $N/T$ is isomorphic to the constant group scheme associated with $W$, i.\,e.
$N/T(R)$ can be identified with a certain subset of the group algebra $RW$, containing $W$,
and $N/T(R)=W$ if there are no nontrivial idempotents in $R$.
An element $\dot w\in N(R)$ is called a representative of $w\in W$
if its canonical image in $N/T(R)$ equals $w$. Obviously, two representatives of an element $w\in W$
differ by an element from $T(R)$.
Therefore, the Gauss cell $\mf G_w=wTUU^-$ does not depend on the choice of a representative
of $w$. A given element $a\in\mf G_w(R)$ can be written as $a=\dot wuv$ for some
$\dot w\in N(R)$, $u\in U(R)$, and $v\in U^-(R)$.

For a root $\alpha\in\Phi$ and an invertible element $\eps$ of $R$ put
$$w_{\alpha}(\eps)=x_\alpha(\eps)x_{-\alpha}(-\eps^{-1})x_\alpha(\eps).$$
It is easy to verify that $w_\alpha(\eps)\in N(R)$ and its image $s_\alpha$ in the Weyl group
is a reflection with respect to the root $\alpha$.
Hence any element $w=\prod_i s_{\alpha_i}\in W$ has a preimage
$\dot w=\prod w_{\alpha_i}(1)\in N(R)$ that comes from $E(\Z)$.
It is well known that it acts on the root elements as follows:
\begin{equation}\label{WeylAction}
x_\alpha(r)^{\dot w}=x_{w(\alpha)}(\pm r)
\end{equation}

For $w\in W$ set
$$
U_w=\langle X_\alpha(R)\mid \alpha\in\Phi^+,\ w(\alpha)\in\Phi^-\rangle
$$
Then $B(R)wB(R)=U(R)wU_w(R)$, and $a=b'\dot w b''$ is called the reduced Bruhat decomposition
of $a\in G(R)$, where
$b'\in U(R)$, $b''\in U_w(R)$, and $\dot w\in N(R)$ is a representative of $w$
(here the element of the Weyl group $w=T(R)\dot w$ is considered as a coset).

It is also known that the big Bruhat cell $Bw_0B$, where $w_0$ is the longest element of
the Weyl group, is a principle open subscheme of $G$, see e.\,g.~\cite[p.\,160]{Jantzen}.
Since $w_0^2=1$ and $w_0Bw_0=B^-$, the Gauss cells
$$
\mathcal G_w=BB^-w=Bw_0B(w_0w)
$$
are shifts of the big Bruhat cells. Hence they are principal open subschemes of $G$ as well.

If $F$ is a field, then the group $G(F)$ splits into the disjoint union of the Bruhat cells $B(F)wB(F)$
over all $w\in W$, see~\cite[2.11]{BorelTits}.
Given $w\in W$ the Bruhat cell is contained in a Gauss cell:
$$
BwB=BwU_w=BU_w^{w^{-1}}w\subseteq BU^-w=\mathcal G_w.
$$
Therefore, over a field $F$ Gauss cells $\mathcal G_w(F)$ cover the group of points $G(F)$.
This means that the family
$$
\{\mathcal G_w\mid w\in W\}
$$
is an open cover of $G$ by principal open subsets.

\smallskip
For a parabolic subgroup $P$ denote by $U_P$ its unipotent radical and by
$L_P$ its Levi subgroup.
Recall that by default all parabolic subgroups are assumed to contain the fixed split maximal
torus $T$. On the other hand, we do not require them to contain the same Borel subgroup.
If an order on the root system is given, then the standard Borel subgroup is the Borel
subgroup, containing all positive root subgroups.
A parabolic subgroup is called standard if it contains the standard Borel subgroup.
Let $\alpha$ be a simple root. Denote by $P_\alpha$ the standard maximal
parabolic subgroup that does not contain the root subgroup $X_{-\alpha}$.

\section{Basic reductions}
In this section we recall some well known facts that reduce the proof of the normal structure theorem
to extraction of a nontrivial element from a parabolic subgroup.
Their proofs with different levels of generality can be found
in~\cite{SteinChevalley,AbeSuzuki,VasersteinChevalley,AbeNormal,VavStav}.

The first assertion is a level computation of a subgroup $H\le G(R)$ normalized by $E(R)$.
In general a level is an admissible pair $(\mf q^s,\mf q^l)$ of additive subgroups of
$R$ such that the intersection of $H$ with the root subgroup $X_\alpha(R)$ equals $X_\alpha(\q^s)$
for a short root $\alpha\in\Phi$ and $X_\alpha(\q^l)$ for a long $\alpha$.
If the structure constants are invertible the situation is much simplier. Namely,
$\mf q^s=\mf q^l$ is an ideal of $R$. To formulate this result put
$\q_\alpha(H)=\{t\in R\,|\,x_\alpha(t)\in H\}$, where $\alpha\in\Phi$.

\begin{lem}\label{level}
Suppose that a root system $\Phi$ and a ring $R$ satisfy Condition~\ref{StrConstants}.
Let $H$ be a subgroup of $G(R)$ normalized by $E(R)$. Then
$\q_\alpha(H)$ is an ideal of $R$ and does not depend on the root $\alpha\in\Phi$.
\end{lem}

In the rest of the article $\q_\alpha(H)$ is denoted by $\q(H)$. It is called the level of the subgroup~$H$.
Thus, for the proof of the normal structure theorem it suffices to show that
$H\le C\bigl(R,\q(H)\bigr)$. Factor out the ideal $\q(H)$. If the image of $H$ lies in the center
we are done. Otherwise, we shall show that $H$ contains a nontrivial root unipotent element.
Using the standard commutator formula
$$
[E(R),G(R,\q)]=E(R,\q)
$$
one lifts this element to a root unipotent element in
$H\smallsetminus E(R,\q)$, which contradicts the definition of the level.

\begin{lem}\label{reduction}
Suppose that for an arbitrary quotient ring $\bar R$ of $R$ any noncentral subgroup 
$H\le G(\bar R)$ normalized by $E(\bar R)$
contains a nontrivial root unipotent element.
Then the normal structure of the group $G(R)$ is standard.
\end{lem}

In the general linear group one extracts a unipotent element as follows. First,
a matrix was converted to a matrix with a zero entry, then one produces a matrix with a ``zero''
column (a column that coincides with the corresponding column of the identity matrix),
and finally an elementary transvection was obtained, see e.\,g.~\cite{GolubchikNormal,BV84}.
An analog of a matrix with a ``zero'' column is an element from a parabolic subgroup.
Indeed, we always can produce a root unipotent element starting from a nontrivial element
from a parabolic subgroup. This statement was proved in~\cite[Theorem~1]{VavStav}

\begin{lem}\label{InP}
Let $H$ be a subgroup of a Chevalley group $G(R)$, normalized by the elementary subgroup
$E(R)$. Suppose that $H\cap P(R)$ is not contained in the center of $G(R)$
for some proper parabolic subscheme $P$ of $G$.
Then $H$ contains a root unipotent element $x_{\alpha}(r)$, $\alpha\in\Phi$,
$r\in R\smallsetminus\{0\}$.
\end{lem}

\section{Extraction of unipotents}

An analog of a matrix with a zero entry is an element, contained in the product of two
parabolic subgroups that are not opposite. Using Lemma~\ref{InPUQ} it is possible to produce
a root unipotent starting from such an element.
For our proof a slightly different case is required. Namely, Lemma~\ref{InPUQ} produces a root
unipotent stating from an element contained in $U_Q(R)P(R)$, where $P$ and $Q$ are arbitrary
parabolic subgroups that contain a fixed split maximal torus.
To ensure that the extracted root unipotent is nontrivial we need to take care about
elements that commute with a root subgroup. This problem is solved by the next statement,
which is well known by specialists.

\begin{lem}\label{centralizer}
Let $a\in G(R)$ be an element such that $[a,X_\alpha(R)]\subseteq C(R)$.
Then $a$ belongs to the group of $R$-points of a proper parabolic subscheme.
\end{lem}

In the sequel we frequently use the following commutator identity, which can be easily verified
by a straightforward calculation. Let $x,y,z$ be elements of an abstract group $G$.
Then
\begin{equation}\label{xyzz-1}
[x,yz]^{z^{-1}}=(x^{-1})^{z^{-1}}x^y=[z^{-1},x]\cdot[x,y].
\end{equation}

\begin{lem}\label{InPUQ}
Let $P,Q$ be two proper parabolic subgroups containing a fixed split maximal torus.
Suppose that $H$ contains a noncentral element from
$U_Q(R)P(R)$. Then $H$ contains a nontrivial root unipotent element.
\end{lem}

\begin{proof}
Clearly we may assume that $P$ is a maximal parabolic subgroup and $Q$ is a Borel subgroup,
in which case $U_Q=U$.
Since the rank of $\Phi$ is at least 2, there exists a root subgroup $X_\alpha$,
contained in the intersection $U\cap L_P$.

By the condition of the lemma there exist $a\in U(R)$ and $b\in P(R)$ such that $ab\in H$.
Let $U^{(0)}(R)=U(R)$ and $U^{(i+1)}(R)=[U^{(i)}(R),U(R)]$.
Since the group $U(R)$ is nilpotent, there exists $k\in\N$ such that
$U^{(k)}_Q(R)=\{1\}$.
Let $i$ be the largest integer such that $a\in U^{(i)}(R)$.
We proceed by induction on $k-i$.

If $i=k$, then $a=1$ and the statement coincides with Lemma~\ref{InP}.
Otherwise, by formula~\eqref{xyzz-1} for an element $r\in R$ we have
$$
[x_\alpha(r),ab]^{b^{-1}}=[b^{-1},x_\alpha(r)]\cdot[x_\alpha(r),a]\in\bigl(P(R)U^{(i+1)}_Q(R)\bigr)\cap H.
$$
If this element is central for any $r$, then a root unipotent element can be extracted by
Lemmas~\ref{centralizer} and~\ref{InP}. Otherwise, the result follows directly from the induction
hypothesis.
\end{proof}

In particular, if $Q=P^-=B$, then the last lemma shows
that a root unipotent element can be extracted from the main Gauss cell.
This partial case allows us to handle subradical subgroups.

\begin{cor}\label{UnderRad}
Let $J$ be the Jacobson radical of $R$. If $H\cap G(R,J)\not\subseteq C(R)$,
then $H$ contains a nontrivial root unipotent element.
\end{cor}

\begin{proof}
A Gauss cell is a principal open subscheme of $G$.
Therefore, an element lies in this cell if and only if it maps a certain element of the
affine algebra to an invertible element of $R$. Since invertibility does not depend on the
Jacobson radical, an element of $G(R,J)$ lies in the same principal open sets as the identity
element $e_R$. On the other hand, $e_R$ obviously belongs to the main Gauss cell.
It follows that $H$ contains a noncentral element from the main Gauss cell and,
by Lemma~\ref{InPUQ}, contains a nontrivial root unipotent.
\end{proof}

The following calculation immediately implies the normal structure theorem of a Chevalley group over a field.

\begin{lem}\label{InCell}
Suppose that $H$ contains a noncentral element from the Gauss cell $U^-(R)B(R)w$.
Then it contains a nontrivial root unipotent.
\end{lem}

\begin{proof}
Suppose that $a=bc\dot w\in H$ for some $b\in U^-(R)$, $c\in B(R)$, and
a preimage $\dot w$ of $w$ in the group $N(R)$.
Let $\alpha\ne\beta$ be simple roots, $P=P_\alpha^w$, and $Q=P^-_\beta$. Then $H$ contains the
element
\begin{multline*}
h=[x_\alpha(r)^{b^{-1}},a]=
x_\alpha(-r)^{b^{-1}} x_\alpha(r)^{b^{-1}a}=\\
x_\alpha(-r)^{b^{-1}} x_\alpha(r)^{c\dot w}\in
Q(R)U_P(R).
\end{multline*}
If this element is noncentral, then by Lemma~\ref{InPUQ} $H$ contains a nontrivial root unipotent.
Otherwise, if $h$ belongs to the center for any $r\in R$, then the element $h^b=[x_\alpha(r),a^b]$
is central as well. Hence, $a^b\in H$ is a noncentral element, which commutes with the root subgroup 
$X_\alpha(R)$ modulo the center.
By Lemma~\ref{centralizer}, $a^b$ lies in a proper parabolic subgroup. Finally, Lemma~\ref{InP}
shows that the subgroup $H$ contains a nontrivial root unipotent.
\end{proof}


\begin{cor}\label{fields}
Let $R$ be a field, satisfying Condition~\ref{StrConstants}, and let $H$ be a subgroup of $G(R)$,
normalized by $E(R)$. Then either $H$ is contained in the center of $G(R)$ or it contains $E(R)$.
\end{cor}

\begin{proof}
If $R$ is field, then any element of the group $G(R)$ belongs to a Gauss cell.
If $H$ contains a noncentral element, then by Lemma~\ref{InCell}, $H$ contains a nontrivial root unipotent.
Now, Lemma~\ref{level} immediately implies that $H$ contains $E(R)$.
\end{proof}

\section{A new look at the decomposition of unipotents}
Recall that $A=\Z[G]$ denotes the affine algebra of a Chevalley--Demazure group scheme
$G=\opn{G}(\Phi,\blank)$ and $g\in G(A)$ is its generic element.
Since $G$ is connected, $A$ is a domain. Denote by $F$ the field of fractions of $A$.
In the sequel we identify elements of $G(A)$ with their canonical images in $G(F)$.
For any element $w\in W$ the Gauss cell $\mf G_w=UB^-w$ is a principal open subscheme of
$G$ (see~\cite[II, 1.9]{Jantzen}), hence $g\in\mf G_w(F)$.
Fix elements $u\in U(F)$, $b\in B^-(F)$ and $\dot w\in N(\Z)$ such that $g=ub\dot w$.

Let $P$ be a proper parabolic subscheme of $G$.
Let $v\in L_P(A)\cap U(A)$. Consider the element $a=v^{u^{-1}}\in U(F)$.
By the dilation principle (see e.\,g.~\cite[Corollary~5.3]{StepComput})
there exists $v$ such that $a\in U(A)$
(actually, since all factors of $a$ belong to $U(F)$, the dilation principle follows easily from the
Chevalley commutator formula).
Clearly, $a^g=v^{b\dot w}\in P^-(F)^{\dot w}\cap G(A)\le P^-(A)^{\dot w}$.
Thus, $a$ is a good element with respect to $g$. In other words, the set of good elements
contains the group $J_w=\bigl(L_P(A)\cap U(A)\bigr)^{u^{-1}}\bigcap E(A)$ for each $w\in W$.
To prove various propositions via decomposition of unipotents it suffices to show that
the span of all subgroups $J_w$ contains at least one root subgroup for each root length.

We have already mentioned in the introduction that for the goal of the current article it
suffices to show that for any ring $R$ the element $h(a)$ does not vanish in $G(R)$
modulo the center for at least one good element $a$ constructed above and $h\in G(R)$.
Denote by $H_g=\langle g\rangle^{E(A)}$ the smallest subgroup of $G(A)$,
containing $g$ and normalized by $E(A)$.

\begin{lem}\label{generic}
Let $P$ be a parabolic subscheme of $G$ and $w\in W$.
There exists an element $c\in H_g$, which lies in the product $P^-(A)^{\dot w}U(A)$.

Moreover, let $S$ be a subscheme defined by the formula
$$
S(R)=\{h\in G(R)\mid h(c)\in C(R)\}.
$$
Then the group of points $S(R)$ is a proper subset of $G(R)$ for any ring $R\ne\{0\}$.
\end{lem}

\begin{proof}
Denote by $s$ an element of the affine algebra $A$ such that $\mf G_w=\opn{Sp}_\Z A_s$.
Decompose $g$ into a product $g=ub\dot w$, where  $u\in U(A_s)$, $b\in B^-(A_s)$, and $\dot w\in N(\Z)$.
Let $\alpha\in\Phi$ be a positive root such that $X_\alpha\le L_P$.
By the dilation principle~\cite[Corollary~5.3]{StepComput}
there exists a positive integer $k$ such that $a=x_\alpha(s^k)^{u^{-1}}$ belongs to
the group $E(A)\cap U(A_s)=U(A)$. Put
$$
c=[g,a^{-1}]=a^ga^{-1}=x_\alpha(s^k)^{b\dot w}a^{-1}\in P^-(A)^{\dot w}U(A)\cap H_g.
$$

Clearly, if $w(\alpha)=\alpha$, then we can choose the element $\dot w$ that commutes with
$X_{\pm\alpha}$.
Now, define $h\in\mf G_w(\Z)$ as follows: if $w(\alpha)=\alpha$, then
$h=x_{-\alpha}(1)\dot w$, otherwise $h=\dot w$.
By the choice of $s$ its image under the homomorphism $h:A\to\Z$ is invertible in $\Z$, i.\,e. equals $\pm1$.
We may assume that $h(s)=1$ (otherwise we can change $s$ to $-s$). Since $h(g)=h$, we have
$
h=h(ub\dot w)=h(u)h(b)\dot w=ex_{-\alpha}(\xi)\dot w,
$
where $\xi$ equals $0$ or $1$. Anyway, the uniqueness of representation of $h$ as a product
of an element from $U(\Z)$ by an element from $B^-(\Z)$ by $\dot w$ implies that $h(u)=e$.
Thus, we get
\begin{align*}
&h(c)=[x_{-\alpha}(1)\dot w,x_\alpha(1)^{-1}]=[x_{-\alpha}(1),x_\alpha(-1)]\text{ or}\\
&h(c)=[\dot w,x_\alpha(1)^{-1}]=x_{w(\alpha)}(\pm 1)x_\alpha(-1).
\end{align*}
Calculation shows that the image of $h(c)$ in $G(R)$ does not lie in $C(R)$ for any
$R\ne\{0\}$. Hence, the image of $h$ in $G(R)$ does not belong to $S(R)$.
\end{proof}

Actually, for the proof of the normal structure theorem it suffices to prove Lemma~\ref{generic}
for one parabolic subscheme $P$ and one element $w\in W$. For instance, if $P$ contains $B$ and $w$ is
the longest element of the Weyl group, then the element $c$ constructed above already belong to
the parabolic subgroup $P(A)$. The freedom in the choice of
$P$ and $w$ is expected to be useful for a description of subgroups of $\opn{G}(\Phi,R)$,
normalized by $\opn{E}(\Delta,R)$, where $\Delta$ is subsystem of $\Phi$ satisfying certain
conditions (see~\cite{VavStepSurvey}).

\section{Proof of the normal structure theorem}\label{NormalSec}

To finish the proof of the normal structure theorem we need the following statement,
which is well known by specialists.

\begin{lem}\label{HallWitt}
Let $H$ be a subgroup of $G(R)$, normalized by $E(R)$. Suppose that if $\Phi=C_2$,
then $R$ has no residue fields of 2 elements.
Then $\bigl[[H,E(R)],E(R)\bigr]=[H,E(R)]$.
In particular, if $[H,E(R)]\le C(R)$, then $H\le C(R)$.
\end{lem}

\begin{proof}
By~\cite[Corollary~4.4]{SteinChevalley} under the assumption of the lemma the group $E(R)$ is perfect,
i.\,e. coincides with derived subgroup. Using the Hall--Witt identity we get
$$
\bigl[E(R),H\bigr]=\bigl[[E(R),E(R)],H\bigr]\le\bigl[[H,E(R)],E(R)\bigr].
$$
The inverse inclusion is obvious. The second assertion follows immediately from the first one and
the fact that the centralizer of $E(R)$ in $G(R)$ coincides with $C(R)$.
\end{proof}

\begin{thm}
Let $G$ be a Chevalley--Demazure group scheme with a reduced irreducible root system $\Phi$ of rank at least 2
and let $R$ be a commutative ring with a unit.
Suppose that $\Phi$ and $R$ satisfy Condition~\ref{StrConstants}. Then, given a subgroup $H\le G(R)$,
normalized by $E(R)$, there exists a unique ideal $\mf a$ in $R$ such that
$$E(R,\mf a)\le H\le C(R,\mf a).$$
\end{thm}

\begin{proof}
By Lemma~\ref{reduction} it suffices to prove that if $H$ is not inside the center, then
it contains a nontrivial root unipotent element. Let $P$ be a proper parabolic subscheme of $G$, $w\in W$, and
$c\in G(A)$ the element from Lemma~\ref{generic}.
If there exists $h\in H$ such that $h(c)\notin C(R)$, then $h(c)\in H$ is a noncentral element from
$P^-(R)^wU(R)$. Then, by Lemma~\ref{InPUQ} $H$ contains a nontrivial root unipotent element.

Otherwise $H$ is contained in the set of $R$-points of the subscheme $S$ defined in Lemma~\ref{generic}.
Consider a maximal ideal $\mf m$ of $R$. The image $\bar H$ of the subgroup $H$ under the canonical
homomorphism $G(R)\to G(R/\mf m)$ is contained in $S(R/\mf m)$ and is normalized by
$E(R/\mf m)$. Since $E(R/\mf m)$ is not contained in $S(R/\mf m)$,
$\bar H$ is contained in $C(R/\mf m)$ by Lemma~\ref{fields}. It follows that the image of the subgroup
$[H,E(R)]$ vanishes in $G(R/\mf m)$, i.\,e. $[H,E(R)]\subseteq G(R,\mf m)$.

Since $\mf m$ is an arbitrary maximal ideal, we have $[H,E(R)]\subseteq G(R,J)$, where $J$ is the Jacobson
radical of $R$. On the other hand, by Lemma~\ref{HallWitt} $[H,E(R)]$ is not contained in the center of $G(R)$.
Hence, by Lemma~\ref{UnderRad} this subgroup contains a nontrivial root unipotent.
This completes the proof.
\end{proof}

\providecommand{\MRhref}[2]{%
  \href{http://www.ams.org/mathscinet-getitem?mr=#1}{#2}
}
\providecommand{\href}[2]{#2}

\end{document}